\documentclass[12pt,reqno]{article}

\usepackage[top=2.5cm,bottom=2.5cm,left=2.5cm,right=2.5cm]{geometry}
\usepackage{graphicx}
\usepackage{latexsym}
\usepackage{amsmath,amssymb}
\usepackage{amscd}
\usepackage[arrow,matrix]{xy}
\usepackage{graphicx}
\usepackage{xcolor}
\usepackage{times}
\usepackage[]{subfigure}
 \usepackage{cite}
 \usepackage{float}

\usepackage{hyperref}

\usepackage[colorinlistoftodos]{todonotes}

\newcommand\func{\operatorname}
\newcommand\grad{\func{grad}}

\usepackage{amscd,amsmath,amsthm,amssymb}

\theoremstyle{plain}
\newtheorem{theorem}{Theorem}[section]

\theoremstyle{definition}
\newtheorem{remark}[theorem]{Remark}

\begin{document}

\title{On timelike Bonnet surfaces in Lorentzian 3-manifold}
\author {{Athoumane Niang$^{1}$}\thanks{{
 E--mail: \texttt{athoumane.niang@ucad.edu.sn} (A. Niang)}},\texttt{ }Ameth  Ndiaye$^{2}$\footnote{{
 E--mail: \texttt{ameth1.ndiaye@ucad.edu.sn} (A. Ndiaye)}}, \texttt{ }Adama Thiandoum$^{3}$\footnote{{
 E--mail: \texttt{adama1.thiandoum@ucad.edu.sn} (A. Thiandoum)}}, \\
\begin{small}{$^{1}$D\'epartement de Math\'ematiques et Informatique, FST, Universit\'e Cheikh Anta Diop, Dakar, Senegal.}\end{small}\\ 
\begin{small}{$^{2}$D\'epartement de Math\'ematiques, FASTEF, Universit\'e Cheikh Anta Diop, Dakar, Senegal.}\end{small}\\
\begin{small}{$^{3}$D\'epartement de Math\'ematiques et Informatique, FST, Universit\'e Cheikh Anta Diop, Dakar, Senegal.}\end{small}}
\date{}
\maketitle%


\begin{abstract} 
In this paper we generalized a result of Soley Ersoy and Kemal Eren \cite{Ersoy} about Bonnet timelike surface in Minkowski 3-space. We give a necessary and sufficient condition for a surface $M$ in a Lorentzian 3-space to be timelike Bonnet surface. At the end, a theorem of classification of timelike Bonnet surface in a Lorentzian 3-space is given.
\end{abstract}
\begin{small} {\textbf{MSC:} 53A10, 53C42, 53C50.}
\end{small}\\
\begin{small} {\textbf{Keywords:} Mean curvature, Principal curvature, Bonnet surface, Lorentzian space.} 
\end{small}\\
\maketitle

\section{Introduction}
The study of submanifolds of a given ambiant space is a naturel interesting problem which enriches our knowledge and understanding of
the geometry of the space itself. For more works about surfaces in three dimensional Lorentzian manifold, see \cite{Calvaruso, Athou1, Athou2, Athou3}.\\
The Bonnet surfaces are surfaces which admit a family of one parameter of  isometric deformations preserving the principal curvatures. V. Lalan (1949) was the first who used the term
“\textit{Bonnet surface}” . In \cite{Bonnet}, Bonnet  proved that any surface with constant mean curvature in Euclidean 3-space $\mathbb{E}$ (which is not totally umbilical) is a Bonnet surface. In \cite{Cartan}, Cartan gave some detailed results for Bonnet surfaces. He proved that in the Euclidean 3-space $\mathbb{E}$, any Bonnet surface is a \textit{Weingarten surface}. If the ambient space of the surface is a complete simply connected Riemannian 3-manifold $\mathbb{R}^3(c)$ of constant curvature $c\geq 0$, W. Chen and H. Li (1997) showed that there exist always Bonnet surfaces which are not Weingarten surfaces \cite{Chen} and in \cite{Chen1}  they obtain the classification results of timelike Bonnet surfaces in 3-dimensional Minkowski space.\\
The results of Bonnet are generalized by Lawson in \cite{Lawson} to any surface with constant mean curvature in Riemannian 3-manifold of constant curvature. In \cite{Chern}, Chern proved that the characterization for isometric deformation preserving the principal curvatures of surfaces is obtained by the aid of differential forms.\\
In \cite{Roussos1, Roussos2, Roussos3} I. Roussos has confirmed that helical surfaces may belong to a class of the Bonnet surface, and he considered the tangent developpabe surfaces as Bonnet surfaces. By using the method of Chern he also obtained a characterization for isometric deformation preserving the mean curvature. More detailed results concerned with these  surfaces are presented \cite{ChenX}.\\
Motivated by the above works, in this papier we give some necessary conditions for a surface to be timelike Bonnet surface in a Lorentzian 3-space $\mathbb{L}^3$. We give a generalization of the works in \cite{Ersoy}.\\
The paper is organized as follow: in Section 2, we recall some basic result and formulas about timelike surface without umbilic points in a three dimensional Lorentzian manifold. In Section three, we study geometric properties of the timelike Bonnet surfaces in $\mathbb{L}^3$.

\section{Some basic formulas}
Let $M$ be a local timelike surface of three dimensional Lorentzian manifold $\mathbb{L}^3$ defined by an embedding $x:M\to\mathbb{L}^3$. The curvature tensor of $\mathbb{L}^3$ is denote by $R$. The metric tensor of $\mathbb{L}^3$ is $g$. We will assume that the surface does not admit umbilic points and has a diagonalizable Weingarten operator $A$. \\
Let $(e_1, e_2)$ be an orthonormal basis on the local surface $M$ which are eigenvector of the shape operator $A$ at any point of $M$. We will assume that $e_1$ is timelike and $e_2$ is spacelike. We put $e_3=e_1\times e_2$, where $\times$ is the vector product of $\mathbb{L}^3$.\\
Let $(\omega^1, \omega^2, \omega^3)$ be the dual basis of the frame $(e_1, e_2, e_3)$ on the surface $M$. The Levi-Civita connection of $\mathbb{L}^3$ and the induced connection on $M$ are both denote by $\nabla$. We have the following equations:
\begin{eqnarray}\label{eq1}
    {\left\{ \begin{array}{cccc}
dx &=& \omega^1e_1+\omega^2e_2  \\
\nabla e_1 &=& \omega_1^2e_2+\omega_1^3e_3\\
\nabla e_2 &=& \omega_2^1e_1+\omega_2^3e_3\\
\nabla e_3 &=& \omega_1^3e_1+\omega_2^3e_2
\end{array}\right.},
\end{eqnarray}
with $\omega^3_3=\omega_1^1=\omega_2^2=0$,  $\omega_1^2=\omega_2^1$, $\omega_1^3=-\omega_3^1$ and $\omega_2^3=\omega_3^2$.\\
We define the shape operator $A:T_pM\to T_pM$ by $A(e_1)=ae_1$, $A(e_2)=ce_2$; so we have 
\begin{eqnarray}\label{eq2}
    \omega_1^3=a\omega^1; \,\,\,\,\,\, \omega_2^3=-c\omega^2.
\end{eqnarray}
The functions $a$ and $c$ are the principal curvature $(a\neq c)$. The mean curvature $H$ and the Gauss curvature $K$ of the surface $M$ are given by $$H=\frac{a+c}{2},$$ $$K=ac$$ respectively and we put $J=\frac{a-c}{2}>0$, for more detail see \cite{Bang, Chen1, Ersoy}.\\
We put 
$$\omega_1^2=h\omega^1+k\omega^2,$$
where $h$ and $k$ are uniquely determined by the first structural equations of $M$ given by 
\begin{eqnarray}\label{eq3}
d\omega^1 = \omega^2\wedge\omega_2^1; \,\,\,\,\,\,d\omega^2 = \omega^1\wedge\omega_1^2. 
\end{eqnarray}
The Gauss and Codazzi equation of the surface $M$ are given respectively by
\begin{eqnarray}\label{eq4}
d\omega_1^2 = -K\omega^1\wedge\omega^2+\psi_1^2,
\end{eqnarray}
with $\psi_1^2(X, Y)=g(R(X,Y)e_1, e_2)$ for all $X, Y\in TM$ and 
\begin{eqnarray}\label{eq5}
    {\left\{ \begin{array}{cc}
d\omega_1^3 - \omega_1^2\wedge\omega_2^3=\psi_1^3  \\
d\omega_2^3 - \omega_2^1\wedge\omega_1^3=\psi_2^3 
\end{array}\right.},
\end{eqnarray}
with $\psi_i^3(X, Y)=g(R(X,Y)e_i, e_3)$, $i=1,2$,  for all $X, Y\in TM$. We will put, 
\begin{eqnarray}\label{eq6}
    \psi_1^3=(c-a)\lambda_1^3\omega^2\wedge\omega^1, \,\,\,\,\,\,\,\psi_2^3=(c-a)\lambda_2^3\omega^1\wedge\omega^2,
\end{eqnarray}
Using the equations (\ref{eq3})-(\ref{eq6}) and the fact that $\omega_1^2=h\omega^1+k\omega^2$ we have 
$$\left[da+(a-c)(h+\lambda_1^3)\omega^2\right]\wedge\omega^1=0.$$ So there exist a function $p'$ such that 
$$da+(a-c)(h+\lambda_1^3)\omega^2=p'\omega^1,$$
and we put $p'=p(c-a)$. So we have
$$da+(a-c)(h+\lambda_1^3)\omega^2=p(c-a)\omega^1.$$
A same computation shows that there exist a function $q$ such that
$$dc+(c-a)(k+\lambda_2^3)\omega^1=(a-c)q\omega^2 .$$
So in summary we get
\begin{eqnarray*}
    {\left\{ \begin{array}{cc}
da+(a-c)(h+\lambda_1^2)\omega^2=(c-a)p\omega^1 \\
dc+(c-a)(k+\lambda_2^3)\omega^1=(a-c)q\omega^2 
\end{array}\right.}.
\end{eqnarray*}
For more information and details see \cite{Ersoy, Chen}.\\
Thus we have these following relations
\begin{eqnarray}\label{eq7}
    {\left\{ \begin{array}{cc}
\frac{da}{c-a}=p\omega^1+(h+\lambda_1^3)\omega^2 \\
\frac{dc}{a-c}=q\omega^2+(k+\lambda_2^3)\omega^1
\end{array}\right.}, 
\end{eqnarray}
and
\begin{eqnarray}\label{eq8}
    2dH=(a-c)\left[ \left((k+\lambda_2^3)-p\right)\omega^1+\left(q-(h+\lambda_1^3)\right)\omega^2\right].
\end{eqnarray}
Now if we put
\begin{eqnarray}\label{eq9}
    {\left\{ \begin{array}{cc}
u= (k+\lambda_2^3)-p\\
v= q-(h+\lambda_1^3)
\end{array}\right.}, 
\end{eqnarray}
the equation (\ref{eq8}) becomes 
\begin{eqnarray}\label{eq10}
    2dH=(a-c)(u\omega^1+v\omega^2).
\end{eqnarray}
We introduce the following 1-forms on the surface $M$
\begin{eqnarray}\label{eq11}
    {\left\{ \begin{array}{cc}
\theta^1=u\omega^1+v\omega^2\\
\theta^2=v\omega^1+u\omega^2
\end{array}\right.}.
\end{eqnarray}
The first fundamental form $I$ of the surface $M$ is given by 
\begin{eqnarray}\label{eq12}
    I=(\omega^1)^2-(\omega^2)^2
\end{eqnarray}
and then (by (\ref{eq10}) and (\ref{eq12})), the gradient of the mean curvature $H$ satisfied
\begin{eqnarray}\label{eq13}
    2\grad H=(a-c)(ue_1-ve_2).
\end{eqnarray}
We introduce again the following 1-forms
\begin{eqnarray}\label{eq14}
    {\left\{ \begin{array}{cc}
\alpha^1=u\omega^1-v\omega^2\\
\alpha^2=-v\omega^1+u\omega^2
\end{array}\right.}.
\end{eqnarray}
We define the $\star$ Hodge operator  by
\begin{eqnarray}\label{15}
    \star\omega^1=\omega^2, \,\,\,\,\,\,\,\star\omega^2=\omega^1, \,\,\,\,(\star)^2=1,
\end{eqnarray}
and the connection form $\omega_1^2$ and the forms $\theta^i$, $\alpha^i$, $i=1,2$ in (\ref{eq11}) and (\ref{eq14}) satisfy
\begin{eqnarray}\label{eq16}
    \star \omega_1^2=k\omega^1+h\omega^2,
\end{eqnarray}
\begin{eqnarray}\label{eq17}
    {\left\{ \begin{array}{cc}
\star\theta^1=\theta^2\\
\star\theta^2=\theta^1
\end{array}\right.}
\end{eqnarray}
and 
\begin{eqnarray}\label{eq18}
    {\left\{ \begin{array}{cc}
\star\alpha^1=\alpha^2\\
\star\alpha^2=\alpha^1
\end{array}\right.}.
\end{eqnarray}
From (\ref{eq7}) and (\ref{eq9}) we have 
\begin{eqnarray}\label{eq19}
    {\left\{ \begin{array}{cc}
\frac{da}{a-c}=\left(u-(k+\lambda_2^3)\right)\omega^1-(h+\lambda_1^3)\omega^2 \\
\frac{dc}{a-c}=(k+\lambda_2^3)\omega^1+\left(v+(h+\lambda_1^3)\right)\omega^2
\end{array}\right.}.
\end{eqnarray}
By using (\ref{eq19}) and (\ref{eq16}), we get 
\begin{eqnarray}\label{eq20}
    dln\vert a-c\vert=\alpha^1-2\star(\omega_2^1+\beta), 
\end{eqnarray}
where $\beta$ statisfies
\begin{eqnarray}\label{eq21}
    {\left\{ \begin{array}{cc}
\beta=\lambda_1^3\omega^1+\lambda_2^3\omega^2\\
\star\beta=\lambda_2^3\omega^1+\lambda_1^3\omega^2
\end{array}\right.}.
\end{eqnarray}
\section{Timelike Bonnet surface in $\mathbb{L}^3$ }
Here, we will study an analogue of the criterion given by \cite{Chen1, Ersoy} for a timelike surface in the Lorentzian manifold $(\mathbb{L}, g)$ to be a Bonnet surface.\\
Let $\Bar{M}$ be a new timelike surface in $\mathbb{L}^3$ with principal curvature $\Bar{a}$ and $\Bar{c}$, which is an isometric deformation of $M$ that preserve the first fundamental form and principal curvature $a$ and $c$ of $M$ (i.e an isometric $\varphi: M\to \bar{M}$ such that $\bar{a}\circ\varphi=a$ and $\bar{c}\circ\varphi=c$). Thus we have
\begin{eqnarray}\label{eq22}
    \Bar{a}=a, \,\,\,\,\,\,\,\,\Bar{c}=c.
\end{eqnarray}
Let $\lbrace \Bar{\omega}^1, \Bar{\omega}^2\rbrace$ be coframe of the frame $\lbrace \Bar{e}_1, \Bar{e}_2\rbrace$ of $\Bar{M}$. Since $\Bar{M}$ and $M$ have the same first fundamental form , we get
\begin{eqnarray}\label{eq23}
    (\Bar{\omega}^1)^2-(\Bar{\omega}^2)^2=(\omega^1)^2-(\omega^2)^2.
\end{eqnarray}
Then there exists some function $\varphi$ on $M$ such that
\begin{eqnarray}\label{eq24}
    {\left\{ \begin{array}{cc}
\Bar{\omega}^1=\cosh{\varphi} \omega^1+\sinh{\varphi} \omega^2\\
\Bar{\omega}^2=\sinh{\varphi} \omega^1+\cosh{\varphi }\omega^2
\end{array}\right.}.
\end{eqnarray}
By taking the exterior differentiation of (\ref{eq24}) with the first structural equation, we can see easily that
\begin{eqnarray*}
    \Bar{\omega}_1^2=\Bar{\omega}_2^1=\omega_1^2-d\varphi.
\end{eqnarray*}
The conditions $\Bar{a}=a$, $\Bar{c}=c$, the relations in (\ref{eq16}), (\ref{eq21}), with (\ref{eq20}) give
\begin{eqnarray}\label{eq25}
    \alpha^1-2\star(\omega_2^1+\beta)=\Bar{\alpha}^1-2\star(\Bar{\omega}_2^1+\Bar{\beta}).
\end{eqnarray}
Applying the operator $\star$ in the equation (\ref{eq25}) with equations (\ref{eq14}-\ref{15}) we get 
\begin{eqnarray}\label{eq26}
    \alpha^2-2(\omega_2^1+\beta)=\Bar{\alpha}^2-2(\Bar{\omega}_2^1+\Bar{\beta}).
\end{eqnarray}
Since $\omega_1^2-\Bar{\omega}_2^1=d\varphi$ and by the equation (\ref{eq25}), the equation (\ref{eq26}) becomes
\begin{eqnarray}\label{eq27}
    \alpha^2-\Bar{\alpha}^2=2d\varphi+2(\beta-\Bar{\beta}).
\end{eqnarray}
The equation $-2dH=(\Bar{a}-\Bar{c})\Bar{\theta}^1$ show that $\Bar{\theta}^1=\theta^1$. So $\Bar{u}\Bar{\omega}^1+\Bar{v}\Bar{\omega}^2=u\omega^1+v\omega^2$; and combining with (\ref{eq24}) give us 
\begin{eqnarray}\label{eq28}
    {\left\{ \begin{array}{cc}
\Bar{u}=u\cosh{\varphi}-v\sinh{\varphi}\\
\Bar{v}=-u\sinh{\varphi}+v\cosh{\varphi }
\end{array}\right.}.
\end{eqnarray}
Now by using the equations (\ref{eq14}) and (\ref{eq28}) with the equation (\ref{eq24}), an easy computation shows that 
\begin{eqnarray}\label{eq29}
    \Bar{\alpha}^2=\sinh{2\varphi}\alpha^1+\cosh{2\varphi}\alpha^2.
\end{eqnarray}
We define $T=\coth{\varphi}$ and we use the relation in (\ref{eq27}) to obtain 
\begin{eqnarray}\label{eq30}
    dT=T\alpha^1+\alpha^2+\frac{\beta-\Bar{\beta}}{\sinh^2{\varphi}}.
\end{eqnarray}
We can rewrite $\lambda_1^3$ and $\lambda_2^3$ given in (\ref{eq6}) by the following expression
\begin{eqnarray}\label{eq31}
    {\left\{ \begin{array}{cc}
\lambda_1^3=\frac{2}{J}g(R(e_1, e_2)e_1, e_3)\\
\lambda_2^3=-\frac{2}{J}g(R(e_1, e_2)e_2, e_3)
\end{array}\right.}.
\end{eqnarray}
By using the equation (\ref{eq24}), we find 
\begin{eqnarray}\label{eq32}
    {\left\{ \begin{array}{cc}
\Bar{e}_1=\cosh{\varphi}e_1-\sinh{\varphi}e_2\\
\Bar{e}_2=-\sinh{\varphi}e_1+\cosh{\varphi}e_2
\end{array}\right.}, 
\end{eqnarray}
and then we get 
\begin{eqnarray}\label{eq33}
    \Bar{e}_1\times\Bar{e}_2=e_1\times e_2=\Bar{e}_3=e_3.
\end{eqnarray}
Using the equations (\ref{eq31}), (\ref{eq32}) and (\ref{eq33}) we get 
\begin{eqnarray}\label{eq34}
    {\left\{ \begin{array}{cc}
\bar{\lambda}_1^3=\left(\cosh{\varphi}\lambda_1^3+\sinh{\varphi}\lambda_2^3\right)\\
\bar{\lambda}_2^3=\left(\sinh{\varphi}\lambda_1^3+\cosh{\varphi}\lambda_2^3\right)
\end{array}\right.},
\end{eqnarray}
and 
\begin{eqnarray}\label{35}
    \Bar{\beta}=\cosh2\varphi \beta+\sinh2\varphi \star\beta.
\end{eqnarray}
One can see that $\beta$ and $\Bar{\beta}$ are related by
\begin{eqnarray}
\frac{\beta-\Bar{\beta}}{\sinh^2\varphi}=-2(\beta+T(\star\beta)).
\end{eqnarray}
Now the equation (\ref{eq30}) becomes 
\begin{eqnarray}\label{eq36}
     dT=T(\alpha^1-2\star\beta)+(\alpha^2-2\beta).
\end{eqnarray}
Let now 
\begin{eqnarray}\label{*}
    {\left\{ \begin{array}{cc}
\gamma^1=\alpha^1-2\star\beta\\
\gamma^2=\alpha^2-2\beta
\end{array}\right.},
\end{eqnarray}
so the equation (\ref{eq36}) become
\begin{eqnarray}\label{**}
     dT=T\gamma^1+\gamma^2.
\end{eqnarray}
\begin{remark}\label{R1}
    The equation (\ref{**}) is a totally differential equation that the hyperbolic anlge $T$ satisfies during the isometric deformation. In order to have non-trivial isometric deformation it is necessary and sufficient that the equation (\ref{**}) to be completely integrable.
\end{remark}
\begin{remark}
When $H$ is constant then $u=v=0$ and $\alpha^1=\alpha^2=0$. Moreover if $M$ has zero normal curvature (i.e $\lambda_1^3=\lambda_2^3=0\,\,\Leftrightarrow\,\, \beta=\star\beta=0$), then $T=constant$ is solution of (\ref{**}).
\end{remark}
So we have the following theorems:
\begin{theorem}
    Let $M$ be timelike surface in $\mathbb{L}^3$ of constant mean curvature $H$ with $H^2>K$. If $M$ has zero normal curvature then $M$ has one parameter family of non trivial isometric deformations preserving the mean curvature; that is $M$ is timelike Bonnet surface in $\mathbb{L}^3$.
\end{theorem}

\begin{theorem}
Let $M$ be timelike surface in $\mathbb{L}^3$. If $M$ has non-constant mean curvature $H$ or has non-zero normal curvature then the exterior differentiation of (\ref{*}) gives
\begin{eqnarray}
Td\gamma^1+(d\gamma^2-\gamma^1\wedge\gamma^2)=0;
\end{eqnarray}
So, the equation (\ref{**}) is completely integrable iff 
$$d\gamma^1=0\,\,\,\,\,\,; \,\,\,\,\,\,d\gamma^2=\gamma^1\wedge\gamma^2.$$
\end{theorem}

\end{document}